\documentclass[10pt,legalblPaper]{article}
\usepackage{amsmath}
\usepackage{amsthm}
\usepackage[utf8]{inputenc}
\usepackage[american]{babel}

\newtheorem{thm}{Theorem}

\title{Does Ten Have a Friend?}
\author{Jeffrey Ward\\wardjm@clarkson.edu\footnote{This research was supported by the National Science Foundation under grant no. 0353723, and was completed during the Auburn University Research Experience for Undergraduates program in 2007.}}
\date{}
\begin{document}
	\maketitle
	\begin{abstract}
		Any positive integer $n$ other than 10 with abundancy index $9/5$ must be a square with at least $6$ distinct prime factors, the smallest being $5$. Further, at least one of the prime factors must be congruent to 1 modulo 3 and appear with an exponent congruent to 2 modulo 6 in the prime power factorization of $n$.
	\end{abstract}
	\section{The Abundancy Index}
	For a positive integer $n$, the sum of the positive divisors of $n$ is denoted $\sigma (n)$; the ratio $\frac{\sigma (n)}{n}$ is known as the \textit{abundancy ratio} or \textit{abundancy index} of $n$, denoted $I(n)$. A \textit{perfect number} is a positive integer $n$ satisfying $I(n) = 2$. 
	
	Considering the millenia-old interest in perfect numbers and the (at least) centuries-old interest in the ``abundancy'' of positive integers, it is somewhat surprising that study of the abundancy index seems to have flourished only relatively recently; see \cite{erdos}, \cite{laatsch}, and \cite{weiner}, and the references there to earlier work. Interesting questions have been asked and answered: for instance, it is now known (\cite{laatsch} and \cite{weiner}) that both the range of the function $I$ and the complement of that range in the rational numbers are dense in the interval $(1, \infty)$.
	
	Questions about another kind of density remain. Let, for $x>1$, $J(x) = I^{-1}((x,\infty)) = \{n \mid I(n) > x \}$; does the limit
	\begin{equation}
		f(x) = \lim_{N \to \infty} \frac{ | J(x) \cap \{ 1,\ldots,N \} | }{N}
	\end{equation}
	exist?
	
	If so, what can be said about the behavior of the non-increasing function $f$? Is it continuous? Strictly decreasing?
	
	The open question about the abundancy index to be addressed here, stated in the title and explained in the next section, is not so exotic -- in fact, it has a claim to the title of least exotic of the unanswered questions about the function $I$. The proof of the main result will use only the most elementary properties of the abundancy index, which we now state with minimal explanation. Proofs may be found in \cite{laatsch} and \cite{weiner}.
	
	\subsection{Elementary Properties of the Abundancy Index}
	Let $m$ and $n$ be positive integers. In what follows, all primes are positive.
	\begin{enumerate}
		\item $I(n) \geq 1$ with equality only if $n=1$.
		\item If $m$ divides $n$ then $I(m) \leq I(n)$ with equality only if $m = n$.
		\item If $p_1,\ldots,p_k$ are distinct primes and $e_1,\ldots,e_k$ are positive integers then
		\begin{align*}
			I( \prod_{j=1}^{k}p_{j}^{e_{j}} ) &= \prod_{j=1}^k\left( \sum_{i=0}^{e_j}p_j^{-i} \right)\\
			&= \prod_{j=1}^k \frac{p_j^{e_j+1}-1}{p_j^{e_j}(p_j-1)}
		\end{align*}
		These formulae follow from their well-known analogues for $\sigma$:
		\begin{align*}
			\sigma(\prod_{j=1}^{k}p_{j}^{e_{j}}) = \prod_{j=1}^{k}\left( \sum_{i=0}^{e_j}p_j^{i} \right) = \prod_{j=1}^k \frac{p_j^{e_j+1}-1}{p_j-1}
		\end{align*}
		Property 3 directly implies a property of $I$ shared by $\sigma$.
		\item $I$ is weakly multiplicative (meaning, if $m$ and $n$ are relatively prime, then $I(mn) = I(m)I(n)$).
		\item Suppose that $p_1,\ldots,p_k$ are distinct primes, $q_1,\ldots,q_k$ are distinct primes, $e_1,\ldots,e_k$ are positive integers and $p_j \leq q_j$, $j = 1, \ldots, k$. Then
		\begin{align*}
			I\left( \prod_{j=1}^k p_j^{e_j} \right) \geq I\left( \prod_{j=1}^k q_j^{e_j} \right)
		\end{align*}
		with equality only if $p_j = q_j$, $j = 1,\ldots, k$. This follows from 3 and the observation that if $e \geq 1$, then $\frac{x^{e+1}-1}{x^e(x-1)}$ is a decreasing function of $x$ on $(1, \infty)$.
		\item If the distinct prime factors of $n$ are $p_1,\ldots,p_k$, then $I(n) < \prod_{j=1}^k \frac{p_j}{p_j-1}$. Although related to $5$, $7$ is most easily seen by applying $3$ and the observation that for $p>1$,
		\begin{align*}
			\frac{p^{e+1}-1}{p^{e+1}-p^e} = \frac{p - \frac{1}{p^e}}{p-1}
		\end{align*}
		increases to $\frac{p}{p-1}$ as $e \rightarrow \infty$.
	\end{enumerate}
	
	\section{Friends}
	Positive integers $m$ and $n$ are \textit{friends} if and only if $m \neq n$ and $I(m) = I(n)$. Thus the perfect numbers form a coterie of friends. This should not be confused with \textit{amicable numbers} -- just because a number is amicable, that doesn't mean it is friendly. Two numbers $m$ and $n$ are \textit{amicable} if and only if $m \neq n$ and $\sigma(m) - m = \sigma(n) - n$. Thus, two numbers are amicable if and only if the sum of their proper divisors are equal. We mention amicable numbers only for clarity; the rest of the article shall be on friendliness.
	
	The ``friends'' terminology seems to have been introduced in \cite{anderson}, where it is asked if the density of the friendly integers -- i.e., positive integers that have at least one friend -- is one. That is, is it true that
	\begin{align*}
		\lim_{N \to \infty} \frac{ | \{ n \mid n~\mathrm{ is ~friendly} \} \cap \{ 1,\ldots,N \} | }{N} = 1?
	\end{align*}
	
	If $m$ and $n$ are friends, and $k$ is a positive integer relatively prime to both $m$ and $n$, then, by the weak multiplicativity of $I$ (property 4), $mk$ and $nk$ are friends. It follows that the friendly integers are fairly numerous. In fact, it is easy to show that for any friendly integer, the set of its friendly multiples has positive (lower) density in the positive integers. Yet the original question in \cite{anderson} remains open.
	
	On the dark subject of unfriendliness, the most elementary observation is that, as a consequence of property 2 of the abundancy index, if $m$ divides $n$, then $m$ and $n$ cannot be friends. Less elementary, but still quite easy to see, is the fact that no prime power has a friend. Therefore, of the numbers $1,2,\ldots,9$, only $6$, a perfect number, has a friend.
	
	Which brings us to the question of the title of this article, asked in \cite{anderson} and in \cite{ford} and is still unanswered. If $10$ does have a friend, the following may be of use in finding it.
	\begin{thm}
		If $n$ is a friend of 10, then $n$ is a square with at least 6 distinct prime factors, the smallest being 5. Further, at least one of $n$'s prime factors must be congruent to 1 modulo 3, and appear in the prime power factorization of $n$ to a power congruent to 2 modulo 6. If there is only one such prime dividing $n$, then it appears to a power congruent to 8 modulo 18 in the factorization of $n$.
	\end{thm}
	\begin{proof}
		Since $I(n) = \frac{\sigma(n)}{n} = I(10) = \frac{9}{5}$, $5\sigma(n) = 9n$ and we see that $5\mid n$. Therefore 2 does not divide $n$, for if it did, 10 would divide $n$, cancelling the possibility of friendship.
		
		Since $n$ is odd and $5\sigma(n) = 9n$, $\sigma(n)$ is odd. As noted in \cite{weiner}, if both $n$ and $\sigma(n)$ are odd, then $n$ must be a square. (To see this, apply the formula $\sigma(\prod_{j=1}^k p_j^{e_j}) = \prod_{j=1}^k \left( \sum_{i=0}^{e_j} p_j^i \right)$ noted after 3 in section 1. If all the $p_j$ are odd and the product is odd then all the $e_j$ must be even.)
		
		If $3 \mid n$ then $n = 3^{2a}\cdot 5^{2b} \cdot m^2$ for positive integers $a,b$, and $m$, with $m$ divisible by neither $2,3$, nor $5$. It is straightforward to check that $I(3^4\cdot 5^2)$, $I(3^2 \cdot 5^4) > \frac{9}{5}$, so by property 2 in Section 1, the only possibility is $a = b = 1$. Then
		\begin{align*}
			9n = 3^4 5^2 m^2 = 5\sigma(n)\\
			= 5\sigma(3^2)\sigma(5^2)\sigma(m^2)\\
			= 5\cdot 13 \cdot 31 \sigma(m^2)
		\end{align*}
		Therefore $13,31 \mid m$, so
		\begin{align*}
			I(n) \ge I(3^2 5^2 13^2 31^2) > \frac{9}{5},
		\end{align*}
		as is easily checked. Therefore, $3$ does not divide $n$.
		
		So $n = 5^{2a}\prod_{i=1}^k p_i^{2e_i}$ for positive integers $a,e_1,\ldots,e_k$, $k \ge 1,$ and distinct primes $p_1,\ldots, p_k > 5$. It is easy to see that $k \ge 4$, for, if $k \leq 3$, then, applying 6 and 7 of Section 1, we would have $I(n) \leq I(5^{2a}7^{2e_1}11^{2e_2}13^{2e_3}) < \frac{5}{4}\frac{7}{6}\frac{11}{10}\frac{13}{12}$, and it is straightforward to check that $\frac{5}{4}\frac{7}{6}\frac{11}{10}\frac{13}{12} < \frac{9}{5}$.
		
		The demonstration that $k \ge 5$ will use 5, 6, and 7 from Section 1.
		
		Verify that $I(5^{2}7^{2}11^{2}13^{2}19^{2}) > \frac{9}{5}$, which implies that $I(5^{2}7^{2}11^{2}13^{2}17^{2}) > \frac{9}{5}$. Verify also that $I(5^{4}7^{2}11^{2}13^{2}23^{2}) > \frac{9}{5}$. Thus if $n = 5^{2a}7^{2e_1}11^{2e_2}13^{2e_3}23^{2e_4}$, then $a = 1$. But then $9n = 5\sigma(n) = 5\cdot 31 \sigma(m^2)$ would imply that $31 \mid n$, which does not hold.
		
		So the cases $k = 4$, $p_1 = 7, p_2 = 11, p_3 = 13$, and $p_4 \in \{17,19,23\}$, are ruled out, and now we lean heavily on 6 and 7 of Section 1 to see that in all other cases when $k = 4, I(n) < \frac{9}{5}$. Thanks to 6, and the cases ruled out so far, only 2 onerous verifications need be performed:
		\begin{align*}
			I(5^{2a}7^{2e_1}11^{2e_2}13^{2e_3}29^{2e_4}) < \frac{5}{4}\frac{7}{6}\frac{11}{10}\frac{13}{12}\frac{29}{28} < \frac{9}{5}
		\end{align*}
		and
		\begin{align*}
			I(5^{2a}7^{2e_1}11^{2e_2}17^{2e_3}19^{2e_4}) < \frac{5}{4}\frac{7}{6}\frac{11}{10}\frac{17}{16}\frac{19}{18} < \frac{9}{5}.
		\end{align*}
		(In each case, the first inequality follows from 7.) Thus, $k \ge 5$.
		
		Finally, since
		\begin{align*}
			5\sigma(n) = 5(1 + \cdots + 5^{2a}) \prod_{i=1}^k \left( \sum_{j=0}^{2e_i} p_i^j \right) = 9n,
		\end{align*}
		we have that $9 \mid \sigma(n)$. If $p \equiv 2 \mod 3$, then $1 + p + \cdots + p^{2e} \equiv 1 \mod 3$ for any positive integer $e$. Consequently, some $p_i \equiv 1 \mod 3$, and $1 + \cdots + p_i^{2e_i} \equiv 0 \mod 3$ implies $2e_i + 1 \equiv 0 \mod 3$. Thus $e_i = 3t +1$ for some integer $t$, so $2e_i = 6t + 2$.
		
		If $p_i$ is the only such prime dividing $n$, then $1 + p_i + \cdots + p_i^{2e_i} \equiv 0 \mod 9$. Checking the possibilities $p_i \equiv 1,4$, or $7 \mod 9$, one finds that $2e_i \equiv 8 \mod 18$.
	\end{proof}
	
	The method of proof of the theorem can, of course, be exploited to get further results, too numerous to mention. For instance, with $k=5$, the method shows that there are only finitely many possibilities for $n$ to check -- and then it's on to $k=6$, unless a friend of 10 has been found with $k=5$.
	
	But the search is long! Intriguingly, if we relax our definitions we can take a shortcut to friendship by going to infinity: observe that $\lim_{k\to \infty} I(3^k5) = \frac{3}{2}\frac{6}{5} = \frac{9}{5} = I(10)$. Of course, since the range of $I$ is dense in $(1,\infty)$, for any positive integer $m$ there are loads of sequences $(n_k)$ such that $\lim_{k\to \infty} I(n_k) = I(m)$, but an inspection of the proofs of that density shows that, generally, the easiest way to come by such sequences is to take the $n_k$ to be products of blocks of very large primes. Let us define a \textit{theoretical friend of proximity} $t$ of a positive integer $m$ to be a sequence $(n_k)$ of positive integers such that $\lim_{k \to \infty} I(n_k) = I(m)$ and $| \{ p \mid p ~\mathrm{is~ a ~positive ~prime ~and,~ for ~some ~} k \mathrm{, ~} p \mid n_k \} | = t$.
	
	Thus, $(3^k5)$ is a theoretical friend of 10 of proximity 2. Does every positive integer have a theoretical friend of finite proximity?\\
	
Acknowledgements: to Mary Drennen, Tracy Gunter, Pete Johnson, Matthew Schneider, and Matthew Ward for comments, observations, and encouragement in the writing of this paper.

\nocite{*}
\bibliography{refs}
\bibliographystyle{plain}

\end{document}